\documentclass[12pt]{article}

\usepackage{amsfonts,latexsym,amssymb,hyperref}

\begin{document}

\newtheorem{defi}{\bf{Definition}}[section]
\newtheorem{lemma}{\bf{Lemma}}[section]
\newtheorem{prop}{\bf{Proposition}}[section]
\newtheorem{theorem}{\bf{Theorem}}[section]
\newtheorem{remark}{\sc{Remark}}[section]

\title{Partial  regularity for the Navier-Stokes-Fourier system}
\author{Luisa Consiglieri\footnote{Independent Researcher Professor, Portugal.
\href{http://sites.google.com/site/luisaconsiglieri}{http://sites.google.com/site/luisaconsiglieri}
}}
\date{}
\maketitle

\begin{abstract}
This paper addresses a nonstationary  flow of heat-conductive
incompressible Newtonian fluid with temperature-dependent viscosity
coupled with linear heat transfer with advection and a viscous heat source
term, under Navier/Dirichlet boundary conditions.
The partial regularity for the velocity of the fluid is proved to each proper
weak solution, that is, for such weak solutions which satisfy some
local energy estimates in a similar way to the suitable weak solutions of 
the Navier-Stokes system. Finally, we study
the nature of the set of
points in space and time upon which proper weak solutions could be singular.
\end{abstract}



\section{Introduction}

Due to diverse applications,
coupled systems governing incompressible flows are subject
of intensive analytical and numerical investigation 
(see for instance \cite{cannon,c09,ljft,diaz,fan,pl1,mrt,zhoufan2} 
and the references therein).
Here we deal with
a nonstationary  flow of heat-conductive
incompressible Newtonian fluid with temperature-dependent viscosity
coupled with linear heat transfer with advection and a viscous heat source
term, under Navier/Dirichlet boundary conditions,
and we study the partial regularity for proper weak solutions to the
coupled system under study.
The proper weak solution is each weak solution that satisfies some
local energy estimates in a similar way to the suitable weak solutions of 
the Navier-Stokes system \cite{proper}.
Although the techniques used in the present work could
be considered standard, the result is new
because the mentioned techniques can not be directly applied.
Under isothermal effects, the viscosity is constant and the problem is described by the Navier-Stokes equations.
We refer to \cite{sch76,sch77,sch80} 
where the partial regularity theory is studied:
at the first instant of time when a viscous incompressible
fluid flow with finite kinetic energy in three space becomes
singular, the singularities in space are concentrated on a closed
set whose one dimensional Hausdorff measure is finite.
This investigation
characterizes some geometric properties and
measures theoretic properties of the sets of
points in space and time upon which weak solutions could be singular 
\cite{ckn,lin}. 
In such conditions,
 it is known that the concept of weak solutions in the sense given by Leray and Hopf
is not sufficient to establish their partial regularity.
 Different regularity criteria for suitable weak solutions to the 
N-S system has been introduced in terms of the smallness of
functionals that are invariant with respect to the natural scaling
either the velocity or its gradient \cite{choe,gk,maha}. 
 We prove the existence of regular points in the sense due to \cite{LS},
 that is, the function is H\"older continuous on a parabolic cylinder
 centred at such point. We remark that in the popular definition
 given at \cite{ckn} the H\"older space is replaced by the space of essentially
 bounded functions.

We recall that the following continuous inclusion
\[W^{2,1}_q(Q_T):=\{v\in L^q(0,T;W^{2,q}(\Omega)):\ \partial_t v \in L^q(Q_T)\}
\hookrightarrow C^{k,\alpha}(\bar Q_T)\]
 only occurs if $q>(n+2)/(2-k)$.
This means for $k=0$ that $q>n/2+1\geq 2$ ($n=2,3$), i.e.,
the Banach space $W^{2,1}_2(Q_T)$ is not embedded in the Banach space of H\"older continuous
functions with exponent $\alpha$ in the $x$-variables and $\alpha/2$ in the $t$-variable.
It is known that the boundedness of the velocity of the fluid, 
 when a  weak solution of N-S equations is care of,
 requires the so-called Ladyzhenskaya-Prodi-Serrin sufficient condition
of the norm of the mixed Lebesgue space $L^{s,r}(Q_T)$ 
 with $2/r+n/s=1$ if $s>n$. Likewise recently 
 a supplementary criterion concerning the N-S equations ensures the
 boundedness of solutions \cite{bae},
 representing the pressure by a relation on the velocity which involves
 its spatial derivatives of  second order.
Indeed,  to guarantee no singularity formation
(so called regularity criteria), there are many
criterion results added on the fluid velocity of weak solutions.
We emphasize that even in 2D, the regularity theory
for the second derivative in space or the derivative 
with respect to the time variable is not available to N-S-F system
because the viscosity is a non-H\"olderian function. It remains to know
whether or not it is possible to find a $L^\infty$-bound
of the velocity of heat-conducting fluids.

 The partial regularity up to the boundary for the N-S-F system
 is still another open problem, since
the study for N-S equations is provided by a constant viscosity
in order to request that the spatial
derivative of second order, 
$\nabla^2 { u }$, belongs to a convenient mixed Lebesgue space,
namely $L^{{9\over 8},{3\over 2}}(Q_T)$ \cite{seregin,tim}
what in the present study does not happen.
Even in the study of parabolic equations 
 the partial regularity up to the boundary is proved by means of the
 existence of $\partial_t u\in L^2(Q_T)$
\cite{ark}. Once more such regularity does not occur in the present study.

Let  $\Omega \subset \mathbb R^n$ be a bounded open domain,  $\partial\Omega$
its boundary, and $T>0$.
Let us consider the initial boundary value problem of the N-S-F system:
\begin{eqnarray}
\begin{array}c
\partial_t {\bf u}- {\rm div }(\mu(\theta) D{\bf u})+({\bf u}\cdot\nabla ){\bf u} ={\bf  f} -
\nabla p\quad\mbox{ in }
Q_T:=\Omega\times]0,T[;\\
{\rm div }~ {\bf u} = 0\quad\mbox{ in }
Q_T;
\end{array}
\label{U1}
\\
\partial_t \theta-k \Delta \theta+{\bf u}\cdot
\nabla\theta =\mu(\theta)  |D{\bf u}|^2\quad\mbox{ in }Q_T;
\label{U2}
\\
{\bf u}\big\vert_{t=0}= {\bf u}_0,\qquad
\theta\big\vert_{t=0}= \theta_0\quad \mbox{ in }\Omega;\label{U3a}\\
u_N:={\bf u}\cdot{\bf n}=0,\quad\tau_T+h(\theta){\bf u}_T={\bf 0},
\quad
\theta= 0\quad\mbox{ on  }\Sigma_T:=\partial \Omega\times ]0,T[,
\label{U3}
\end{eqnarray}
where 
$\theta$ represents   the temperature, $\bf u$  the
velocity of the fluid and
$D{\bf u}={1\over2}(\nabla {\bf u}
 +\nabla{\bf u}^T)$, $p$ denotes the pressure,  $\mu$ the viscosity,
$\bf f$ denotes the given external body forces, $k$ denotes
the conductivity assumed to be a fixed positive constant,
 and ${\bf u}_0$  and $\theta_0$ are some given functions.
For simplicity,
the density and the heat capacity are constants assumed equal to one
and we do not consider the existence of the external heat
source, since the heating dissipative term is the main mathematical difficulty.
  The product of two tensors is given by $D:\tau=D_{ij}\tau_{ij}$
  (in Einstein's convention) and the norm by $|D|^2=D:D$.
 The
boundary conditions in (\ref{U3}) deal with
 the linear Navier law describing the slip fluid-boundary interaction
and for the sake of clarity we found convenient that 
is assumed homogeneous Dirichlet condition for the temperature.
 Here  ${\bf n}=(n_i)$
denotes the unit outward normal to $\partial\Omega$,
$u_N,{\bf u}_T$ are the normal and the tangential components
of the velocity vector, respectively, $\tau_T=\tau\cdot {\bf
n}-\tau_N{\bf n}$ is the tangential component of the deviator
 stress tensor $\tau=\mu(\theta) D{\bf u}$, and
$h$ denotes the friction coefficient. 

The organization of the paper is as follows. Next section is concerned to the
presentation of the appropriate functional framework and the main results.
Section \ref{s3} is devoted to 
the proof of some auxiliary results introduced at Section \ref{ar}.
In Section \ref{r4}, we prove  the partial regularity result (Theorem \ref{main}).
Finally an estimate to the parabolic Hausdorff dimension on the
set of the singularities for the fluid velocity is provided in Section \ref{fim}.
The optimal estimate is still an open problem.

\section{Assumptions and main results}

Here we assume that $\Omega \subset
 \mathbb R^n$ ($n=2,3$) is a bounded open domain  sufficiently regular, e.g. of class
 $ C^2$.
For $1\leq q\leq \infty$, we introduce 
\begin{eqnarray*}
{\bf V}_{q}(\Omega)=\{{\bf v}\in{\bf W}^{1,q}(\Omega):\ \nabla\cdot{\bf v}=0
\mbox{ in }\Omega,\
v_N=0\mbox{ on  }\partial\Omega\};\\
{\bf W}^{1,q}_N(\Omega)=\{{\bf v}\in{\bf W}^{1,q}(\Omega):\ 
v_N=0\mbox{ on  }\partial\Omega\}
\end{eqnarray*}
with norm $$\Vert\cdot\Vert_{1,q,\Omega}=\Vert\nabla\cdot\Vert_{q,\Omega}+
\Vert\cdot\Vert_{q,\partial\Omega}.$$
We denote by bold the vector spaces of vector-valued or tensor-valued functions. 
For any set $A$, we  write $(u,v)_A:=\int_A uv$ whenever
$u\in L^q(A)$ and $v\in L^{q'}(A)$, where $q'=q/(q-1)$ is the conjugate exponent to $q$,
or simply $(\cdot,\cdot)$ whenever there exists no confusing  at all, and we use
 the symbol
$\langle\cdot,\cdot\rangle$ to denote
 a generic duality pairing,
not distinguished between scalar and vector fields.
 
\begin{defi}\label{defnsf}
We say that the triple $({\bf u},p,\theta)$ is a weak solution to the Navier-Stokes-Fourier (N-S-F) problem 
 (\ref{U1})-(\ref{U3}) in $Q_T$ if,  $p\in L^{(n+2)/n}(Q_T)$,
\begin{eqnarray*}
{\bf u} \in {\mathcal U}:=L^{\infty}(0,T;{\bf L}^2(\Omega))\cap L^2(0,T;
{\bf V}_2(\Omega)), \\
\theta \in  {\mathcal E}:=L^{\infty}(0,T;L^1(\Omega))\cap L^q(0,T;W^{1,q}_0 (\Omega)),\quad q<{n+2\over n+1},\\
\partial_t{\bf u} \in {\mathcal X}:=L^{2}(0,T;({\bf W}^ {1,2 }_N(\Omega))')\cap 
L^ {(n+2)/n }(0,T;({\bf W}^{1,(n+2)/2}_N(\Omega))'),\\
\partial_t\theta \in  L^1(0,T;W^{-1,\ell} (\Omega)),\quad {1\over \ell}={n\over q( n+1)}+{n\over 2(n+2)},
\end{eqnarray*}
and satisfies the variational formulation
\begin{eqnarray}
\langle \partial_t {\bf u},{\bf v}\rangle+\int _{Q_T} \Big(\mu(\theta
)D{\bf u}:D{\bf v}+({\bf u}\cdot\nabla ) {\bf u}\cdot{\bf v} \Big) dxdt+\nonumber\\
+\int_{\Sigma_T}h(\theta){\bf u}_T\cdot{\bf v}_T \; dSdt =
\int _{Q_T} ( {\bf f}\cdot {\bf v}+ p\ {\rm div~}{\bf v})dxdt,\label{wmotion}\\
\forall {\bf v} \in L^\infty(0,T;{\bf W}^{1,\infty}_N(\Omega)),\qquad
{\bf u}\big\vert_{t=0}={\bf u}_0\mbox{ in }\Omega;\nonumber\\
\langle \partial_t \theta,\phi\rangle+\int _{Q_T} \Big( k \nabla \theta-\theta
{\bf  u}\Big)\cdot \nabla \phi dx dt
=\int _{Q_T} \mu(\theta)| D{\bf u}|^2  \phi dx dt,\label{wheat}
\\
\forall \phi \in L^{\infty}(0,T; W^{1,\infty}_0(\Omega)),\qquad
\theta\big\vert_{t=0}= \theta_0 \mbox{ in }\Omega.\nonumber
\end{eqnarray}
\end{defi}

In (\ref{wmotion}) the convective term verifies
$ ({\bf u}\cdot\nabla){\bf u}\in{\bf L}^{(n+2)/(n+1)}({Q_T})$,
for every ${\bf u}\in{\mathcal U}\hookrightarrow {\bf L}^{2(n+2)/n}({Q_T})$.
In (\ref{wheat}) the advection term  $ \theta{\bf u}\in{\bf L}^\ell({Q_T})$ if
$\theta\in {\mathcal E}\hookrightarrow L^ {q(n+1)/n}({Q_T})$ 
for $\ell\geq 1$, i.e. $q>1$ if $n=2$ and $q\geq 15/14 $ if $n=3$
($q\geq 2n(n+2)/((n+4)(n+1))$).
For $n=2,3$,
we remark that the embedding holds
${\mathcal X}\hookrightarrow {L}^{(n+2)/n}(0,T;({\bf W}^{1,(n+2)/2}_N(\Omega))')$.
In conclusion, all terms in the above equalities
are meaningful.

Next we define the existence of weak solutions such that verify some
local energy inequalities (the proof of the existence result
can be found in \cite{proper}, for the two-dimensional case).
\begin{defi}
We say that the N-S-F problem, defined by (\ref{wmotion})-(\ref{wheat}), has proper weak solutions if
the weak solution in accordance to Definition \ref{defnsf}
satisfies the following local energy inequalities 
\begin{eqnarray}
\int _\Omega {| {\bf u}(x,t)-{\bf a}|^2\over 2}\varphi^2(x,t)dx+
\int _{Q_t=\Omega\times ]0,t[} \mu(\theta
)|D{\bf u}|^2 \varphi^2 dx\, d\tau\leq\nonumber\\
\leq 2\int _{Q_t}  p\varphi{\bf u}\cdot\nabla\varphi dx d\tau+
2\int_{Q_t}\mu(\theta)\varphi D{\bf u}:(({\bf u}-{\bf a})\otimes\nabla \varphi) dxd\tau+\nonumber\\
+\int_{Q_t}{|{\bf u}-{\bf a}|^2}
\varphi(\partial_t\varphi+{\bf u} \cdot
\nabla\varphi)dxd\tau
+\int_{Q_t}{\bf f}\cdot( {\bf u}-{\bf a})
\varphi^2 dx d\tau,\qquad
\label{e1}\\
\int _\Omega (\theta\psi)(x,t)dx
\leq\int _ {Q_t} \theta
\left( \partial_t\psi+k\Delta\psi+
{\bf u}\cdot \nabla\psi\right)dxd\tau+\nonumber\\
+\int _{Q_t} \mu(\theta)| D{\bf u}|^2\psi dxd\tau,
\label{e2}
\end{eqnarray}
for all ${\varphi} \in {C}^{\infty}_0(Q_T)$, for any ${\bf a}\in\mathbb R^n$,
for all $\psi \in C^{\infty}_0(Q_T)$ such that $\psi\geq 0$, a.e $t\in ]0,T[$. 
Moreover, $\theta\geq 0$ in $Q_T$.
\end{defi}

\begin{remark}
The local energy inequality for the temperature, (\ref{e2}),
plays an essential role in the proof of the decay lemma
(cf. Theorem \ref{T1_1}) because we obtain 
$\theta \in L^{\infty}(0,T;L^1_{\rm loc}(\Omega))$.
\end{remark}

We assume that the following hypotheses hold
\begin{description}
\item[(A1)] ${\bf f}:Q_T\rightarrow\mathbb R^n$ is given such that ${\bf f}\in {\bf L}^{2}(Q_T)$;
\item[(A2)] 
$\mu ,h: {\mathbb R}\to {\mathbb R}$  are continuous functions such that
\begin{eqnarray}
0<\mu_\# \le \mu(s)\le \mu^\# ,\qquad \forall s\in {\mathbb R};
\label{defmu}\\
\label{defan}
0<h_\#\leq h(s)\leq h^\#,\quad \forall s\in{\mathbb R};
\end{eqnarray}
\item[(A3)]
$ {\bf u}_0\in {\bf L}^{2}(\Omega),$ $ \theta_0\in  L^{1}(\Omega)
$ such that
\begin{equation}\label{defeo}
{\rm div }~ {\bf u}_0 = 0\quad\mbox{ in }\Omega;\quad{\rm ess}\inf_{x\in \Omega}\theta_0(x)\geq 0.
 \end{equation}
\end{description}

In order to establish the partial regularity, let us introduce some
additional notations.
\begin{defi}
Given a positive number $\lambda$, we set
\begin{equation}\label{cl}
c_\lambda(f,\omega):=\sup\{{1\over R^{\lambda-2}}
\left(-\hspace*{-0.45cm}\int
_{Q(z,R)}|f|^2\right)^{1/2}:\quad R>0,\quad Q(z,R)\subset\subset \omega\}
\end{equation}
for all $\omega\subset Q_T$. In particular, we set $c_\lambda(f)=c_\lambda(f,Q_T)$.
\end{defi}

Given a positive number $\lambda$ 
we introduce the "parabolic" variant of the Morrey spaces \cite{LS}
\begin{eqnarray*}
M_{2,\lambda}(Q_T)=\{f\in L^2_{\rm loc}(Q_T):\ 
c_\lambda(f,Q_T)<\infty\}.
\end{eqnarray*}

Let  $z_0=(x_0, t_0)\in Q_T$ and $R>0$ be such that
  $Q(z_0, R)\equiv B(x_0,R)\times ]t_0-R^2,t_0[\subset Q_T$.
For any set $A$, 
 we denote by $|A|$ the usual Lebesgue measure of the set,
and we  write $ -\hspace*{-0.35cm}\int_Av:={1\over |A|}\int_A v$ whenever
$v\in L^1(A)$, 
or simply $(v)_A= -\hspace*{-0.35cm}\int_Av$.
Observing that $2(n+2)/(n+4)<a'<q(n+1)/n$ is equivalent to
\[{q(n+1)\over q(n+1)-n}< a<{2(n+2)\over n}=\left\{
\begin{array}{ll}
4&\mbox{if }n=2\\
10/3&\mbox{if }n=3
\end{array}\right.
\]
we can choose $1<q<(n+2)/(n+1)$ such that $q>3n/(2(n+1))$, i.e.
$q(n+1)/(q(n+1)-n)<3$. So we set
\begin{equation}\label{defa}
3\leq a<2(n+2)/n.
\end{equation}
Let us  introduce the following scaling invariant functional:
$$
\begin{array}c
\bar Y_R (z_0;{\bf u},p,\theta) \equiv  \Big( -\hspace*{-0.35cm}\int
 _{Q(z_0, R)} |{\bf u}|^{a}dz\Big)^{1/a}
  + R\Big( -\hspace*{-0.35cm}\int
 _{Q(z_0, R) }|p|^ {n+2\over n}dz\Big)^{n\over n+2} + \\ +   R\Big( -\hspace*{-0.35cm}\int
 _{Q(z_0, R) } |\theta  |^{a'}~dz\Big)^{1/a'}.
\end{array}
$$
The choice of the above exponents is consequence of
${\bf u}\in {\bf L}^{2(n+2)/n}(Q_T)$,
$p\in {L}^{(n+2)/n}(Q_T)$, $\theta \in L^{q(n+1)/n}(Q_T)$
and ${\bf u}\theta\in {\bf L}^{1}(Q_T)$.
We refer \cite{LS} in order to compare to the
those exponents in the  analysis 
 for the N-S system into the three-dimensional space $(n=3)$.

\begin{theorem}\label{main}
Suppose $\lambda$ and $a$ are positive numbers such that (\ref{defa}) is verified.
Let ${\bf f}\in M_{2,\lambda}(Q_T)$ and (A1)-(A3) be fulfilled.
Assume that
$({\bf u},p,\theta)$
is a proper weak solution of  the N-S-F system in
$Q_T$.
Then,  there exists a positive
constant $\Lambda$ {\bf depending only on}  $\mu_\#$, $\mu^\#$, $\lambda$ and $n$.
Moreover, if for any $z_0\in Q_T$
\begin{equation}\label{limsup}
\limsup_{R \to 0^+}R\bar Y_R (z_0; {\bf u},p,\theta)\leq \Lambda,
\end{equation}
then $z_0$ is a regular point in the following sense: the function
$z\mapsto  {\bf u}(z)$ 
is H\" older continuous in some neighbourhood of the point $z_0$.
\end{theorem}

The local energy inequality for the temperature has a nonstandard format. It is still an open problem the proof that the pointwise criterion
\[ 
\limsup_{R \to 0}
 {1\over R}\left(\int
_{Q(z_0,R)}|\nabla{\bf u}|^2
dz\right)^{1/2}
+\left(\int_{Q(z_0,R)}
 |\nabla\theta|^{q}
dz\right)^{1/q}
\leq \gamma_*
\] 
implies the local condition for the regularity (\ref{limsup}),
for some constant $\gamma_*>0$.

\begin{remark}\label{r1}
Under all conditions of Theorem \ref{main}, we denote by $Q_0$
the set of all regular points such that satisfy (\ref{limsup}).
By definition, $Q_0$ is an open set.
The set $S= Q_T\setminus Q_0$, known as the singular set,  reads
 \[
 S=\{z\in Q_T:\ \limsup_{R \to 0^+}R\bar Y_R (z; {\bf u},p,\theta)>\Lambda\}.
\]
 Moreover,  it can be proved that it is  of Lebesgue measure zero.
\end{remark}
 
The parabolic Hausdorff dimension on the singular set 
can be estimated.
Indeed, we state the following result.
 \begin{theorem}\label{dim}
 Let $a$ be given as  in (\ref{defa}). Then
${\mathcal H}^{n-(a-2)/(a-1)}(S)=0,$
where $\mathcal H^d$ is the $d$-dimensional Hausdorff measure
 with respect to  the standard parabolic metric 
of a set $S$
defined as follows
\[
\mathcal{H}^d(S):=\lim_{\delta\rightarrow 0}\inf_{}\{
\sum_{i=1}^\infty R_i^d:\quad 
S\subset \bigcup_{i=1}^\infty Q(z_i,R_i),
\quad 0<R_i\leq \delta\}.
\]
Therefore, the parabolic Hausdorff dimension of $S$ does not exceed $n-(a-2)/(a-1)$,
i.e. dim$_{\mathcal H}(S):=
\inf\{d\geq 0:\ \mathcal{H}^d(S)=0\}\leq n-(a-2)/(a-1).$
\end{theorem}

\section{Auxiliary results}
\label{ar}

The main tool for the partial regularity analysis is the decay property  of the scaled Lebesgue norms
of the triple velocity-pressure-temperature which is based on a standard "blow up" method and the decomposition
of the pressure.
Let us  introduce the following scaling invariant functional:
$$
\begin{array}c
 Y_R (z_0; {\bf u},p,\theta) \equiv 
 \Big( -\hspace*{-0.35cm}\int
 _{Q(z_0,R) } |{\bf u}(x,t)- ({\bf u})_{Q(z_0,R)}|^{a}dz\Big)^{1/a} +\\ + R \Big(
-\hspace*{-0.35cm}\int _{Q(z_0,R ) } |p(x,t)-(p)_{B(x_0,R)}(t)|^ {(n+2)/n}
dz\Big)^{n/(n+2)} + \\ +
R\Big( -\hspace*{-0.35cm}\int _{Q(z_0,R ) } | \theta(x,t) - (\theta)_{Q(z_0,R)}
|^{a'}dz\Big)^{1/a'}.
\end{array}
$$

\begin{theorem}[Decay estimate]\label{T1_1}
Suppose that $0<\varsigma< 1$ and $0<\beta<\lambda$ are fixed numbers.
For each function $\mu$ satisfying (\ref{defmu}) there exist  positive
constants $\gamma$ and $\varepsilon$ depending only on   $\mu_\#$, $\mu^\#$, $\varsigma$, $\beta$ and $\lambda$ such that
for any proper weak solution $({\bf u}, p, \theta)$
of  the N-S-F problem in $Q_T$ verifying
\begin{equation}\label{dh}
\left.\begin{array}{r}
\forall Q(z,R)\subset\subset Q_T:\quad 0<R<\varepsilon,\\
R|( {\bf u})_{Q(z,R)}|<1,\quad R|( {\theta})_{Q(z,R)}|<\gamma/4\\
Y_R(z; {\bf u},p, \theta)+c_\lambda({\bf f}) R^\beta < \gamma
\end{array}\right\}\end{equation}
we have the decay estimate
$$
Y_{\varsigma R}(z; {\bf u},p, \theta) \le C_*
\varsigma^\alpha (Y_R(z;{\bf u},p,\theta)+c_\lambda({\bf f}) R^\beta)
$$
with $0<\alpha<2-n/2$ and some absolute constant $C_*$ depending only on $\mu_\#$ and
$\mu^\#$.
\end{theorem}

To prove Theorem \ref{main} some iterative estimates are  required 
\cite[Lemmas 2.5, 3.1]{LS}.
Let us state an iterative result, in which a careful choice of $\varsigma$
significantly simplifies the iterative formula as well as its proof
(compare to \cite[Lemmas 2.5, 3.1]{LS}).

\begin{lemma}\label{ite}
Under all conditions of Theorem \ref{T1_1}, let $0<\delta<\alpha.$
If additionally $\varsigma$ is such that
\begin{equation}\label{tau}
C_* \varsigma^{\alpha-\delta}\leq 1/2\qquad \mbox{and}\qquad
\varsigma^\delta+\varsigma^\beta\leq 1, 
\end{equation}
then every proper weak solution $({\bf u}, p, \theta)$
of  the N-S-F system in $Q_T$ satisfying (\ref{dh}) verifies
\begin{enumerate}
\item for all $m\in\mathbb N$,
\begin{equation}\label{ind}
Y_{\varsigma^m R}(z; {\bf u},p, \theta) \le \varsigma
^{m\delta} (Y_R(z;{\bf u},p,\theta)+c_\lambda({\bf f}) R^\beta);
\end{equation}
\item for all $\rho\in ]0,\varsigma R]$,
\begin{equation}\label{irho}
Y_{\rho}(z; {\bf u},p, \theta) \le C\left({\rho\over R}\right)^{\delta} (Y_R(z;{\bf u},p,\theta)+c_\lambda({\bf f}) R^\beta),
\end{equation}
where $C=C(\varsigma,\delta)$ denotes a positive constant.
\end{enumerate}
\end{lemma}

Finally, we state the additional technical result.
\begin{lemma}\label{T2}
Assume that $({\bf u},p,\theta)$
satisfies the system (\ref{U1})-(\ref{U2})
 in the sense of distributions.
For every $R,r>0$,  $e_0=(y_0,s_0)$ and  $z_0=(x_0, t_0)$ such that $ Q(z_0, R)\subset Q_T$,
\begin{enumerate}

\item 
if we introduce the functions
$$
\begin{array}c
{\bf u}^R(y, s) \equiv {R\over r} {\bf u}(x_0+{R\over r}(y-y_0), t_0+\left({R\over r}\right)^2( s-s_0)), \\
p^R(y, s) \equiv \left({R\over r}\right)^2 p(x_0+{R\over r}(y-y_0), t_0+\left({R\over r}\right)^2( s-s_0)), \\
\theta^R(y, s) \equiv \left({R\over r}\right)^2 \theta(x_0+{R\over r}(y-y_0), t_0+\left({R\over r}\right)^2( s-s_0)),
\end{array}
$$
then $({\bf u}^R,p^R, \theta^R)$ satisfy the transported system in $Q(e_0,r)$ 
\begin{eqnarray}
\partial_s {\bf u}^R- {\rm div }_y(\mu^R(\theta^R) D_y{\bf u}^R)+({\bf u}^R\cdot\nabla_y
 ){\bf u}^R ={\bf  f}^R -
\nabla_y p^R;
\label{ns}\\
{\rm div }_y {\bf u}^R = 0;\\
\partial_s \theta^R- k\Delta_y  \theta^R+{\bf u}^R\cdot
\nabla_y\theta^R =\mu^R(\theta^R)  |D_y{\bf u}^R|^2; \label{heat}
\end{eqnarray}
with the function $\mu^R(\vartheta)\equiv
\mu(r^2\vartheta /R^2) $ 
 and 
\[
{\bf f}^R(y,s)=\left({R\over r}\right)^3{\bf f}(x_0+{R\over r}(y-y_0),t_0+\left({R\over r}\right)
^2( s-s_0));
\]
\item 
then
\begin{equation}
r\bar Y_r( e_0; {\bf u}^R,p^R, \theta^R)= R\bar Y_R( z_0; {\bf u},p, \theta).\label{IF}
\end{equation}
Moreover
\begin{equation}\label{cr}
c_\lambda({\bf f}^R,Q(e_0,r))\leq \left({R\over r}\right)^{\lambda+1}c_\lambda({\bf f},Q(z_0,R)).
\end{equation}
\end{enumerate}

\end{lemma}
\begin{remark} The
function $\mu^R$ satisfies (\ref{defmu}) with the same constants
$\mu_\#$, $\mu^\#$. 
\end{remark}

Henceforth 
the symbol $C$ will denote a positive, finite
constant that may vary from line to line; the relevant dependencies
 on the data will be specified whenever it will be required.
Notice that the dependence never occurs
 on the unknown functions $\bf u$, $p$ or $\theta$.

\section{Proof of the auxiliary results}
\label{s3}

\subsection{Proof of Lemma \ref{T2}}

1)
Considering the change of variables
\begin{eqnarray*}
Q(e_0,r)\equiv B(y_0,r)\times(s_0-r^2,s_0)\rightarrow Q(z_0,R)\equiv
B(x_0,R)\times(t_0-R^2,t_0)\\
(y,s)\mapsto \left(x_0+{R\over r}(y-y_0), t_0+\left({R\over r}\right)^2 (s-s_0)\right)
\end{eqnarray*}
we get (\ref{ns})-(\ref{heat}), observing that
\begin{eqnarray*}
&&\partial_s {\bf u}^R+({\bf u}^R\cdot\nabla_y){\bf u}^R+\nabla_y\, p^R=\left({R\over r}\right)^3(\partial_t
{\bf u}+({\bf u}\cdot \nabla){\bf u}+\nabla p);\\
&&{\rm div}_y(\mu^R(\theta^R)D_y{\bf u}^R)
=\left({R\over r}\right)^2{\rm div}_y(\mu(\theta)D_y{\bf u})=\left({R\over r}\right)
^3{\rm div}(\mu(\theta)D{\bf u});\\
&&\partial_s {\theta}^R+{\bf u}^R\cdot\nabla_y{\theta}^R-\Delta_y
\theta^R=\left({R\over r}\right)^4(\partial_t
{\theta}+{\bf u}\cdot \nabla{\theta}+\Delta\theta);\\
&&\mu^R(\theta^R)|D_y{\bf u}^R|^2=\left({R\over r}\right)^4\mu(\theta)|D{\bf u}|^2.
\end{eqnarray*}

2) Considering the change of variables
\begin{eqnarray*}
Q(z_0,R)&\rightarrow &Q(e_0,r)\\
(x,t)&\mapsto &\left(y_0+{r\over R}(x-x_0),s_0+\left({r\over R}\right)^2(t-t_0)\right)\end{eqnarray*}
the Jacobian is $(r/R)^{n+2}$. Then, for every $m\geq 1$ and $u\in L^m(Q_T)$,
\[
-\hspace*{-0.45cm}\int_{Q(e_0,r)}|u(x_0+{R\over r}(y-y_0),t_0+\left({R\over r}\right)^2(s-s_0))|^m
de=
  -\hspace*{-0.45cm}\int_{Q(z_0,R)}|u|^m dz.
  \]
In particular, it follows
\begin{eqnarray*}
r({-\hspace*{-0.45cm}\int}_{Q(e_0,r)}|{\bf u}^R|^{a}de)^{1/a}
={ R}({-\hspace*{-0.45cm}\int}_{Q(z_0,R)}|{\bf u}|^{a}dz)^{1/a};\\
r^2(-\hspace*{-0.45cm}\int
_{Q(e_0,r)}|p^R|^{(n+2)/n}de)^{n/(n+2)}=R^2(-\hspace*{-0.45cm}\int_{Q(z_0,R)}|p|^{(n+2)/n
}dz)^{n/(n+2)};\\
r^2(-\hspace*{-0.45cm}\int_{Q(e_0,r)}|\theta^R|^{a'}de)^{1/a'}
=R^2(-\hspace*{-0.45cm}\int_{Q(z_0,R)}|\theta|^{a'}dz)^{1/a'}.
\end{eqnarray*}
Consequently, (\ref{IF}) holds.
Also we have
\[{1\over r^{\lambda-2}}
\left(-\hspace*{-0.45cm}\int_{Q(e_0,r)}|{\bf f}^R|^2
de\right)^{1/2}=\left({R\over r}\right)^{\lambda+1}{1\over R^{\lambda-2}}
  \left(-\hspace*{-0.45cm}\int_{Q(z_0,R)}|{\bf f}|^2 dz\right)^{1/2}.
\]
Thus we obtain (\ref{cr}) 
and Lemma \ref{T2} is proved.

\subsection{Proof of the decay estimate (Theorem \ref{T1_1})}

Suppose the opposite, i.e., there exist a sequence $\{(\varepsilon_m,\gamma_m)\}$
and solutions $({\bf u}_m,p_m,\theta_m)$ 
of (\ref{U1})-(\ref{U2}) in $\Omega_m\times ]0,T_m[$ and $Q(z_m,R_m)\subset\subset 
\Omega_m\times ]0,T_m[$ such that
\begin{eqnarray}\label{cm}
R_m|( {\bf u}_m)_{Q(z_m,R_m)}|<1,\quad R_m|( {\theta}_m)_{Q(z_m,R_m)}|<\gamma_m/4;\\
\label{zero}
\varepsilon_m\rightarrow 0, \quad
Y_{R_m}(z_m;{\bf u}_m,p_m,\theta_m)+d_m
 R_m^\beta= \gamma_m\to 0, \quad (m\rightarrow \infty);\\
Y_{\varsigma R_m}(z_m;{\bf u}_m,p_m,\theta_m)\ge C_* \varsigma^\alpha \gamma_m, \label{Opposite}
\end{eqnarray}
 where $d_m=c_\lambda({\bf f}_m;\Omega_m\times ]0,T_m[)$.
Applying the change of variables
\begin{eqnarray*}
Q\equiv Q(0,1)&\leftrightarrow& Q_m\equiv Q(z_m,R_m),\\
e=(y,s)&\leftrightarrow &z=(x,t)
\end{eqnarray*}
we introduce the transported functions:
\begin{eqnarray*}
{\bf w}_m(e)& :=& \frac 1{\gamma_m} ( {\bf u}_m (z)- ({\bf u}_m)_{Q_m}), \\
\pi_m(e) &:= &{R_m\over \gamma_m} ( p_m (z)- (p_m)_{B_m}),\qquad B_m=B(x_m,R_m), \\
 \varkappa_m(e) &:= &{R_m\over {\gamma_m}} (
\theta_m (z)- (\theta_m)_{Q_m}).
\end{eqnarray*}
These functions satisfy the system (in the sense of distributions) in $Q$
\begin{equation}
\left\{
\begin{array}c
\partial_s {\bf w}_m +R_m\left((\gamma_m {\bf w}_m+{\bf a}_m)\cdot \nabla_y\right){\bf w}_m
- \mbox{\rm div}_y(\mu({ \gamma_m \over R_m}
\varkappa_m +b_m)D_y{\bf w}_m)=\\=-  \nabla_y \pi_m +
{\bf g}_m,\qquad
 {\rm div }_y {\bf w}_m =0, \\
 \partial_s \varkappa_m +R_m(\gamma_m {\bf w}_m+{\bf a}_m)\cdot\nabla_y\varkappa_m-  k\Delta_y\varkappa_m=\\
= R_m\gamma_m  \mu({ \gamma_m\over R_m} \varkappa_m+b_m ) | D_y {\bf w}_m |^2.
\end{array} \right.
\label{NS2}
\end{equation}
Here ${\bf g}_m(e)={R_m^2\over \gamma_m}{\bf f}_m(z) $, ${\bf a}_m:= ({\bf u}_m)_{Q_m}$ and  $b_m:= (\theta_m)_{Q_m}$.

 From (\ref{zero}) we can extract subsequences, still denoted by the same symbols, such that
\[
  {\bf w}_m\rightharpoonup {\bf w} \mbox{ in } {\bf L}^{a}(Q),
\quad \pi_m\rightharpoonup \pi \mbox{ in } L^{(n+2)/n}(Q), \quad
\varkappa_m\rightharpoonup \varkappa \mbox{ in } L^{a'}(Q)
\]
and taking into account that by definition
$({\bf w}_m)_Q=(\varkappa_m)_Q=0$ and $(\pi_m)_B(s) =0,$ $s\in(-1,0)$
and also the mean integrals remain zero for the weak limits,
it follows
\begin{eqnarray}
1\geq \liminf_{m\rightarrow\infty} Y_1(0;{\bf w}_m,\pi_m,\varkappa_m)\geq
Y_1(0;{\bf w},\pi ,\varkappa);\label{y1}\\
\liminf\limits_{m\to\infty} Y_\varsigma(0;{\bf w}_m,\pi_m,\varkappa_m) \ge C_*\varsigma^\alpha.
\label{Contradiction}
\end{eqnarray}

In order to obtain more information about $({\bf w},\pi ,\varkappa)$ we wish to pass to
the limit the system satisfied by $({\bf w}_m,\pi_m,\varkappa_m)$ considering the weak
formulation
\begin{eqnarray}
-\int_Q {\bf w}_m\cdot \partial_s {\bf v}de+\int_Q \mu ({\gamma_m\over R_m}
\varkappa_m+b_m)D_y{\bf w}_m:D_y{\bf v}de=\nonumber\\
=R_m\int_Q(\gamma_m {\bf w}_m+{\bf a}_m) \otimes 
{\bf w}_m:\nabla_y {\bf v}de
+\int_Q \pi_m{\rm div}_y\, {\bf v}de+\int_Q{\bf g}_m
\cdot{\bf v}de ,\nonumber\\ \quad\mbox{ for all } {\bf v}\in{\bf C}_0^\infty(Q);
\label{w1}
\\-\int_Q \varkappa_m \partial_s\phi de-k\int_Q
\varkappa_m\Delta _y\phi de=R_m\int_Q \varkappa_m(\gamma_m {\bf w}_m+{\bf a}_m)\cdot\nabla
_y\phi de+\nonumber\\+R_m  \gamma_m\int_Q
\mu({\gamma_m\over R_m} \varkappa_m+b_m)|D_y{\bf w}_m|^2\phi de,\quad
\mbox{ for all } \phi\in C_0^\infty(Q).\label{w2}
\end{eqnarray}

The dependence of the viscosity on the temperature does not allow to proceed as in the N-S
system (under isothermal behaviour).
 In the following we use the local energy inequalities (\ref{e1})-(\ref{e2})
  for $({\bf u}_m,p_m,\theta_m)$
in order to obtain some indispensable estimates
on the required solutions under some
standard conditions on the functions involved therein. 

\subsubsection{$ L^ {\infty }(-\varsigma^2,0;{\bf L}^{2}(B(0,\varsigma)))\cap
 L^{2}(-\varsigma^2,0;{\bf H}^1(B(0,\varsigma)))$ estimate for ${\bf w}_m$}
The local energy inequality (\ref{e1}) for $({\bf u}_m,p_m,\theta_m)$ can be rewritten for
$({\bf w}_m,\pi_m,\varkappa_m)$, for a.e. $s\in  ]-1,0 [$,
\begin{eqnarray}
\int _B {| {\bf w}_m(y,s)|^2\over 2}\varphi^2(y,s)dy+\int_{-1}^s\int _B\mu({\gamma_m\over R_m}\varkappa_m+b_m
)|D_y{\bf w}_m|^2 \varphi^2 dy\, d\tau\leq\nonumber\\
\leq R_m\int_{-1}^s\int_B \mu({\gamma_m\over R_m}\varkappa_m+b_m)\varphi
D_y{\bf w}_m:({\bf w}_m\otimes\nabla_y \varphi) dyd\tau+\nonumber\\+\int_{-1}^s
\int_B{|{\bf w}_m|^2}\varphi(\partial_s\varphi+R_m(\gamma_m{\bf w}_m+{\bf a }_m) \cdot
\nabla_y\varphi)dyd\tau+\nonumber\\
+2\int _{-1}^s\int_B  \pi_m\varphi{\bf w}_m\cdot\nabla_y\varphi dyd\tau+
\int_{-1}^s\int_B{\bf g}_m\cdot {\bf w}_m\varphi^2dyd\tau,\quad
\forall {\varphi} \in {C}^{\infty}_0(Q),\label{we1}
\end{eqnarray}
where  $B\equiv B(0,1) \subset \mathbb R^n$. 
Using the Cauchy-Schwarz and Young inequalities we deduce
\begin{eqnarray}
\int _B {| {\bf w}_m(y,s)|^2\over 2}\varphi^2(y,s)dy+{1\over 2}\int_{-1}^s\int _B\mu({\gamma_m\over R_m}\varkappa_m+b_m
)|D_y{\bf w}_m|^2 \varphi^2 dy\, d\tau\leq\nonumber\\
\leq {1\over 2}\int_{-1}^s\int_B \mu({\gamma_m\over R_m}
\varkappa_m+b_m)|{\bf w}_m\otimes\nabla_y \varphi|^2dyd\tau+\nonumber\\
 +\int_{-1}^s
\int_B{|{\bf w}_m|^2}\varphi(\partial_s\varphi+R_m(\gamma_m{\bf w}_m +{\bf a}_m)\cdot
\nabla_y\varphi)dyd\tau+\nonumber\\
+2\int _{-1}^s\int_B  \pi_m\varphi{\bf w}_m\cdot\nabla_y\varphi dyd\tau
+\int_{-1}^s\int_B{\bf g}_m\cdot {\bf w}_m\varphi^2
dyd\tau.\label{vw}
\end{eqnarray}
Thus, taking a conveniently chosen $\varphi$ and using (\ref{y1}) under (\ref{defa})
and also (\ref{cm}) we obtain 
\begin{equation}\label{cotawm}
{\rm ess}\sup_{s\in ]-\varsigma^2,0[} \| {\bf w}_m (s)\|_{2, B(0,\varsigma)}+
 \|\nabla_y {\bf w}_m \|_{2, Q(0,\varsigma)} \leq C+ \| {\bf g}_m \|_{2, Q},
\end{equation}
with $C=C(\mu_\#,\mu^\#)$.

\subsubsection{$ L^ {\infty }(-\varsigma^2,0;L^{1}(B(0,\varsigma)))$ 
estimate for ${\gamma_m\over R_m}\varkappa_m+b_m$
and $\varkappa_m$}

The local energy inequality (\ref{e2}) for $({\bf u}_m,p_m,\theta_m)$ can be rewritten for
$({\bf w}_m,\pi_m,\varkappa_m)$, for a.e. $s\in  ]-1,0 [$,
\begin{eqnarray*}
\int _B \left({\gamma_m\over R_m}\varkappa_m+b_m\right)\psi(y,s)dy
\leq \gamma_m \int_{-1}^s\int_B\left(\gamma_m\varkappa_m+R_mb_m\right){\bf w}_m
\cdot \nabla_y\psi dyd\tau +\nonumber\\
+  \int_{-1}^s\int_B\left({\gamma_m\over R_m}\varkappa_m+b_m\right)
\left( \partial_s\psi+k\Delta_y\psi+
{\bf a}_m\cdot \nabla_y\psi
\right)dyd\tau+\\+\gamma_m^2
\int_{-1}^s\int_B \mu({\gamma_m\over R_m}\varkappa_m+b_m)| D_y{\bf w}_m|^2  \psi dyd\tau.
\end{eqnarray*}
Thus, taking a conveniently chosen $\psi$ and using the Gronwall Lemma,
(\ref{cotawm}), (\ref{y1}), (\ref{cm}) and (\ref{defmu}), we obtain 
\begin{equation}\label{cotatm}
{\rm ess}\sup_{s\in ]-\varsigma^2,0[} \|
 {\gamma_m\over R_m}\varkappa_m+b_m \|_{1, B(0,\varsigma)}\leq \gamma_m^2
 C( 1+ \| {\bf g}_m \|_{2, Q})\exp[ C+R_m \gamma_m  \| {\bf g}_m \|_{2, Q}].
\end{equation}
By definition (\ref{cl}) we have
\begin{eqnarray}
 \left(-\hspace*{-0.45cm}\int_{Q}|{\bf g}_m|^2 de\right)^{1/2}
= {R_m^2\over \gamma_m}  \left(-\hspace*{-0.45cm}\int_{Q(z_m,R_m)}|{\bf f}_m|^2 dz\right)^{1/2}
\nonumber\\
\leq {R_m^2\over \gamma_m}d_m R_m^{\lambda-2}
= {R_m^\beta\over \gamma_m}d_m R_m^{\lambda-\beta}
\leq  R_m^{\lambda-\beta}\longrightarrow 0,\quad \mbox{as }m\rightarrow 0\quad (\beta<\lambda).\label{gm}
\end{eqnarray}

Finally, since $b_m\geq 0$ and using (\ref{gm}), from (\ref{cotatm}) we find
\begin{equation}
{\rm ess}\sup_{s\in ]-\tau^2,0[} \| \varkappa_m \|_{1, B(0,\tau)}\leq
   C R_m\gamma_m.\label{vm}
\end{equation}

\subsubsection{Passage to the limit as $m\rightarrow\infty$}

In order to be in conditions to pass to the limit in (\ref{w1})-(\ref{w2}), let us state the following convergences.
From (\ref{cm}) we get $R_m{\bf a}_m\rightarrow {\bf a}\mbox{ in }\mathbb R^n$.
Moreover, $|{\bf a}|<1$ arises.
Inserting (\ref{gm}) into (\ref{cotawm}) we get
$$D_y{\bf w}_m\rightharpoonup D_y{\bf w} \quad \mbox{in }{\bf L}^2(Q(0,\varsigma)), $$
and consequently
\begin{eqnarray*}
 R_m \gamma_m
\int_Q\mu({\gamma_m \over R_m}\varkappa_m+b_m)|D_y{\bf w}_m|^2\phi de
&\leq& R_m\gamma_m\mu^\#\|D_y{\bf w}_m\|^2_{2,Q}\|\phi\|_{\infty,Q }\\
&\leq& CR_m\gamma_m\rightarrow 0.
\end{eqnarray*}
Inserting (\ref{gm}) in (\ref{cotatm})
and (\ref{vm}) as tending $m$ to infinity, we obtain
\begin{equation}\label{bm}
{\gamma_m \over R_m}\varkappa_m+b_m\rightarrow 0,\quad \varkappa_m\rightarrow 0
\qquad \mbox{in }L^1(Q(0,\varsigma))
\end{equation}
 and consequently $\mu({\gamma_m\over  R_m}\varkappa_m+b_m)\rightarrow \mu(0)$
in $\mathbb R$.

Then we can pass to the limit in (\ref{w1})-(\ref{w2}) 
with $\varkappa=0$ and $({\bf w},\pi  )$ is a
solution to the time dependent Stokes system  
(in the sense of distributions)
\[
\left.
\begin{array}c
\partial_s {\bf w}- \mu(0) \Delta_y {\bf w}+\nabla_y \pi =-({\bf a}\cdot\nabla_y){\bf w}=-\nabla_y\cdot({\bf a}\otimes{\bf w})
\\
 {\rm div}_y {\bf w} =0
\end{array} \right\}
\quad\mbox{in }Q. 
\]

The classical theory for the nonstationary Stokes equations \cite{s76} claims
from $({\bf a}\cdot\nabla_y){\bf w }\in {\bf L}^2(Q)$ that
${\bf w }\in {\bf W}^{2,1}_{2}(Q) 
\hookrightarrow L^{2(n+2)/n}(-1,0;{\bf W}^{1,{2(n+2)/n}}(B))$.
Applying the bootstrap argument 
from $({\bf a}\cdot\nabla_y){\bf w }\in {\bf L}^{2(n+2)/n}(Q)$ it follows that
${\bf w }\in{\bf W}^{2,1}_{2(n+2)/n}(Q) \hookrightarrow{\bf C}^{0,\alpha}(\bar Q)$
for $0\leq \alpha <2-n/2$.
Then $\bf w$ is H\"older continuous in the closure of the cylinder $Q(0,\varsigma/2)$
and the following estimate holds
\begin{equation}
Y_\varsigma(0;{\bf w},0 ,0 ) \le C_{**} \varsigma^\alpha
\label{ya}\end{equation}
with some constant $C_{**}$ depending only on ${\mu_\#}$ and $ {\mu^\#}$.

\begin{remark}
The above bootstrap argument can still be applied 
if we use first the embedding
$\mathcal U\hookrightarrow {\bf L}^{2(n+1)/n}(Q)$
and next the $L^p$-theory for the Stokes equation with RHS in the divergence form
\cite{ks}, that is,
from ${\bf w}\in {\bf L}^{2(n+1)/n}(Q) $
it  follows  that ${\bf w}\in
L^ {2(n+1)/n  }(-1,0;{\bf W}^{1,2(n+1)/n}(B))$.
\end{remark}

On the other hand, we need to extract subsequences which converge
strongly in order to pass to the limit in the integral $Y_\varsigma(0;{\bf w}_m,\pi_m, \varkappa_m)$.
From (\ref{w1}) and using (\ref{defmu}), (\ref{cm}), (\ref{cotawm}), (\ref{gm}),
(\ref{y1}), we derive that
$\{\partial_s {\bf w}_m\}_{m\in\mathbb N}$ is bounded in $L^{(n+2)/n}(-1,0;
{\bf W }^{-1,(n+2)/2}(B))$
since the following estimate holds
\begin{eqnarray*}
\int_Q {\bf w}_m\cdot \partial_s {\bf v}de\leq
\left(
\mu^\#\|D_y{\bf w}_m\|_{2,Q}
+R_m\gamma_m \|{\bf w}_m\|_{2(n+2)/n,Q}^2
+\right.\\\left.
+R_m|{\bf a}_m |
\|{\bf w}_m\|_{2(n+2)/n,Q}+\| \pi_m\|_{(n+2)/n,Q}\right)\|\nabla_y {\bf v}\|_{(n+2)/2,Q}
+\|{\bf g}_m\|_{2,Q}\|{\bf v}\|_{2,Q},
\end{eqnarray*}
for every $ {\bf v}\in
L^{(n+2)/2}(-1,0;{\bf W}_0^{1,(n+2)/2}(B))$.
Thanks to (\ref{cotawm}) and (\ref{gm}), the sequence 
$\{{\bf w}_m\}$ is bounded in
$L^{\infty}(-\varsigma^2,0;{\bf L}^2(B(0,\varsigma)))\cap L^2(-\varsigma^2,0;
{\bf H}^{1}(B(0,\varsigma)))
\hookrightarrow{\bf L}^{2(n+2)/n}(Q(0,\varsigma))$.
Using $H^1(B)\hookrightarrow\hookrightarrow L^5(B)$
and a compactness result \cite{si} we obtain
\begin{equation}\label{wm}
{\bf w}_m\rightarrow {\bf w} \quad \mbox{in }{\bf L}^{a}(Q(0,\varsigma)).
\end{equation}

Therefore, using (\ref{wm}), (\ref{ya}) and (\ref{vm})
 we obtain
$$
\begin{array}c
 \limsup\limits_{m\to\infty} Y_\varsigma(0;{\bf w}_m,\pi_m,\varkappa_m)\leq \\
\leq Y_\varsigma
( 0;{\bf w},0,0)+ \limsup\limits_{m\to\infty} \varsigma
 \Big(-\hspace*{-0.35cm}\int
_{Q(0,\varsigma ) } |\pi_m-(\pi_m)_{B(0,\varsigma)}|^{(n+2)/n}de\Big)^{n/(n+2)}\le
\\ \le  C_{**} \varsigma^\alpha + 4\limsup\limits_{m\to\infty} \varsigma
 \Big(-\hspace*{-0.35cm}\int
_{Q(0,\varsigma ) } |\pi_m|^{(n+2)/n}de\Big)^{n/(n+2)}.
\end{array}
$$

In order to estimate the pressure, the idea is to decompose the pressure into
 two parts, which one part is harmonic and hence is smooth (see \cite{LS}).
Choosing ${\bf v}(y,s)=\chi(s)\nabla\varphi(y)$, where $\chi\in C^\infty_0(0,1)$ and $\varphi\in C^\infty_0(B)$,
as a test function in (\ref{w1}), remarking that 
div\,${\bf w}_m=0$ and
\[-\int_Q {\bf w}_m\cdot\nabla\varphi {d\chi\over ds}dy\,ds=\int_Q\nabla\cdot{\bf w}_m\varphi{d\chi\over
ds}dy\,ds=0,\]
we obtain, for a.e. $s\in ]-1,0[$,
\begin{eqnarray*}
\int_B \mu ({\gamma_m\over R_m}
\varkappa_m+b_m)D_y{\bf w}_m:D_y{\nabla_y\varphi}dy-R_m\gamma_m \int_B
{\bf w}_m\otimes 
{\bf w}_m:\nabla_y ^2\varphi dy=\\
=\int_B \pi_m{\Delta}_y {\varphi}dy+\int_B{\bf g}_m
\cdot\nabla_y\varphi dy.
\end{eqnarray*}

Let $q_m$ be defined by
\begin{eqnarray}
\int_B \mu ({\gamma_m\over R_m}
\varkappa_m+b_m)D_y{\bf w}_m:D_y{\nabla_y\varphi}dy-R_m\gamma_m\int_B
 {\bf w}_m \otimes 
{\bf w}_m:\nabla_y ^2\varphi dy=\nonumber\\
=\int_B q_m{\Delta}_y {\varphi}dy+\int_B{\bf g}_m
\cdot\nabla_y\varphi dy,\label{qm1}
\end{eqnarray}
which is valid for a.e. $s\in ]-1,0[$ and for all ${\varphi}\in W^{2,(n+2)/2}_0(B)$,
and
\begin{equation}
(\pi_m-q_m,\Delta _y\varphi)=0.\label{qm2}
\end{equation}
In (\ref{qm2}), if we consider $\varphi$
as a solution of the Laplace problem 
\begin{eqnarray*}
\Delta_y \varphi(s)=|(\pi_m-q_m)(s)|^{2/n}{\rm sign} ((\pi_m-q_m)(s)) &\mbox{in }& B,
\end{eqnarray*}
we obtain $\pi_m\equiv q_m$ a.e. in $Q$.
In (\ref{qm1}), if we consider $\varphi$
as the unique solution of the Dirichlet-Laplace problem 
\begin{eqnarray*}
\Delta_y \varphi(s)=|q_m(s)|^{2/n}{\rm sign} (q_m(s)) &\mbox{in }& B(0,\varsigma);\\
\varphi(s)=0&\mbox{on }&\bar B\setminus B(0,\varsigma),
\end{eqnarray*}
then using (\ref{defmu}) and the H\"older inequality, it results
\begin{eqnarray*}
\int_{B(0,\varsigma)} |q_m|^{(n+2)/n}\leq \|{\bf g}_m
\|_{2,B(0,\varsigma)}\|\nabla_y\varphi\|_{2,B(0,\varsigma)}+\nonumber\\
+ \left(\mu^\#
\|\nabla_y{\bf w}_m\|_{(n+2)/n,B(0,\varsigma)}
+R_m\gamma_m\| {\bf w}_m\|_{2(n+2)/n,B(0,\varsigma)}^2\right)
\|\nabla ^2\varphi\|_{(n+2)/2,B(0,\varsigma)}.
\end{eqnarray*}
It is known that the following estimate holds
\[
\|\nabla_y\varphi(s)\|_{{(n+2)/ 2},B(0,\varsigma)}
+\|\nabla^2_y\varphi(s)\|_{ {(n+2)/ 2},B(0,\varsigma)}
\leq C\|q_m(s)\|_{(n+2)/n,B(0,\varsigma)}^{2/n},
\]
for some constant $C>0$. So it follows
\[
\|\pi_m\|_{{(n+2)/n},B(0,\varsigma)}^{}\leq C\left({R_m}^{\lambda-\beta}+ \mu^\#
\|\nabla_y {\bf w}_m\|_{(n+2)/n,B(0,\varsigma)}+R_m\gamma_m \right),
\]
taking into account (\ref{gm}) and (\ref{zero}).
Applying the H\"older inequality,  we conclude that
\begin{equation}
\limsup\limits_{m\to\infty} -\hspace*{-0.45cm}\int_{Q(0,\varsigma ) }|\pi_m|^{(n+2)/n}de
\leq C_1
\limsup\limits_{m\to\infty}\left( -\hspace*{-0.45cm}\int_{Q(0,\varsigma ) }|
\nabla_y{\bf w}_m|^2de\right)^{n+2\over 2n},
\label{pwm}
\end{equation}
with $C_1=C(\mu^\#)$.

Taking  in (\ref{vw}) $\varphi\in C^\infty_0(Q(0,2\varsigma
))$ the cut-off function such that
$\varphi\equiv 1\mbox{ in }Q(0,\varsigma)$, $|\nabla\varphi|\leq
C/\varsigma ^{2\alpha/(n+2)}$ and $|\partial_t\varphi|\leq C/\varsigma^{4\alpha
/(n+2)}$ in $Q(0,2\varsigma)$, it results
\begin{eqnarray*}
\mu_\#\int_{Q(0,\varsigma)}|\nabla_y{\bf w}_m|^2 dy\, ds 
\leq {C\over\varsigma^{4\alpha/(n+2)}}\int_{Q(0,2\varsigma)}|{\bf w}_m|^2dy\,ds+\\
+{C\over \varsigma^{2\alpha/(n+2)}}\int _{Q(0,2\varsigma)}|\pi_m{\bf w}_m|dy\,ds
+\int_{Q(0,2\varsigma)}|{\bf g}_m|^2dy\,ds.
\end{eqnarray*}
Using (\ref{gm}), the H\"older inequality and observing that
\begin{eqnarray*}
\left({1\over \mu_\#}{C\over\varsigma^{2\alpha/(n+2)}}
-\hspace*{-0.45cm}\int_{Q(0,2\varsigma ) }|\pi_m{\bf w}_m|de\right)^{n+2\over 2n}\leq
{1\over 2C_1}-\hspace*{-0.45cm}\int_{Q(0,2\varsigma ) }|\pi_m|^{(n+2)/n}de+\\
+{C\over\varsigma^{2\alpha/n}}
\left(-\hspace*{-0.45cm}\int_{Q(0,2\varsigma ) }|{\bf w}_m|^{a}de \right)^{n+2\over
an},
\end{eqnarray*}
we conclude
\begin{eqnarray}\label{cotap}
\limsup\limits_{m\to\infty} 
\left(-\hspace*{-0.45cm}\int_{Q(0,\varsigma ) }|\nabla_y{\bf w}_m|^2de \right)^{n+2
\over 2n}\leq 
{C_2\over\varsigma^{2\alpha/n}}
\limsup\limits_{m\to\infty} 
\left(-\hspace*{-0.45cm}\int_{Q(0,2\varsigma ) }|{\bf w}_m|^a de \right)^{n+2
\over an}\nonumber\\
+{1\over 2C_1}
\limsup\limits_{m\to\infty} 
-\hspace*{-0.45cm}\int_{Q(0,2\varsigma ) }|\pi_m|^{n+2\over n}de,
\end{eqnarray}
with $C_2=C(\mu_\#,\mu^\#)$.

Introducing (\ref{cotap}) into (\ref{pwm}) and applying (\ref{ya}),
we obtain 
$$
\limsup\limits_{m\to\infty} 
-\hspace*{-0.45cm}\int_{Q(0,\varsigma ) }|\pi_m|^{(n+2)/n}de
\leq C_{3} (2\varsigma)^\alpha+{1\over 2}
\limsup\limits_{m\to\infty} -\hspace*{-0.45cm}\int_{Q(0,2\varsigma)}|\pi_m|^{(n+2)/n}de,
$$
with $C_3=C_1C_22^{2\alpha/n}C_{**}^{(n+2)/n}$.

Iterating over $i=1,\cdots,k$, where $k$ is such that $1/2<2^k\varsigma<1$ we derive
$$
\limsup\limits_{m\to\infty} -\hspace*{-0.45cm}\int_{Q(0,\varsigma ) }|\pi_m|^{n+2\over
n}de
\leq 4C_{3} \sum_{i=1}^k 2^{i(\alpha-2)}\varsigma^\alpha+
\limsup\limits_{m\to\infty} -\hspace*{-0.45cm}\int_{Q(0,2^k\varsigma ) }|\pi_m|^{n+2
\over n}de
.
$$
Considering that 
$$
 -\hspace*{-0.45cm}\int_{Q(0,2^k\varsigma ) }|\pi_m|^{(n+2)/n}de
\leq C 2^{n+2}\int_Q|\pi_m|^{(n+2)/n}de\leq C 2^{n+2}
$$
it results
$$
\limsup\limits_{m\to\infty}  -\hspace*{-0.45cm}\int_{Q(0,\varsigma ) }
|\pi_m|^{(n+2)/n}de
\leq {4C_{3} \over 2^{2-\alpha}-1}\varsigma^{\alpha}+C
\leq C_{***}^{(n+2)/n}.
$$

 Hence we obtain
$$
\limsup\limits_{m\to \infty} Y_\varsigma (0;{\bf w}_m, \pi_m, \varkappa_m) \leq
C_{**}\varsigma^\alpha+C_{***}\varsigma \leq
(C_{**}+C_{***})\varsigma^\alpha.
$$
This is a contradiction with (\ref{Contradiction}) if we set $C_*
> C_{**}+C_{***}$. Theorem \ref{T1_1}  is proved.

\subsection{Proof of Lemma \ref{ite}}

(i) We prove (\ref{ind}) by induction on $m$.
The proof of  (\ref{ind}) when $m=1$ follows from Theorem \ref{T1_1},
observing that (\ref{dh}) holds and
\[ C_* \varsigma^\alpha=C_* \varsigma^{\alpha-\delta}\varsigma^\delta\leq\varsigma^\delta.\]
Now we suppose that (\ref{dh}) holds and, for 
$i=1,\cdots,m$, we have
\begin{equation}\label{ii}
Y_{\varsigma^i R}(z; {\bf u},p, \theta) \le \varsigma
^{i\delta} (Y_R(z;{\bf u},p,\theta)+c_\lambda ({\bf f})R^\beta).
\end{equation}
Thus (\ref{ii}) and (\ref{dh}) imply that
\begin{eqnarray*}
Y_{\varsigma ^iR}(z; {\bf u},p, \theta) +c_\lambda({\bf f}) (\varsigma^i R)^\beta
&\leq& \varsigma^{i\delta}(Y_R(z;{\bf u},p,\theta)+c_\lambda({\bf f}) R^\beta)
+c_\lambda({\bf f}) \varsigma ^{i\beta}R^\beta\\
&\leq &\varsigma^{i\delta}Y_R(z;{\bf u},p,\theta)+c_\lambda({\bf f}) R^\beta(\varsigma^{i\delta}+\varsigma^{i\beta})<\gamma,
\end{eqnarray*}
considering that (\ref{tau}) holds.

We can apply Theorem \ref{T1_1} and subsequently (\ref{ii}) resulting
\begin{eqnarray*}
Y_{\varsigma ^{i+1}R}(z; {\bf u},p, \theta)& =&Y_{\varsigma(\varsigma
 ^iR)}(z; {\bf u},p, \theta) 
\leq C_*\varsigma^\alpha\left(Y_{\varsigma^iR}(z;{\bf u},p,\theta)+c_\lambda({\bf f}) (\varsigma
^i R)^\beta\right)\\
&\leq &C_*\varsigma^\alpha\left(
\varsigma^{i\delta}(Y_R(z;{\bf u},p,\theta)+c_\lambda ({\bf f})R^\beta)+
c_\lambda({\bf f}) \varsigma ^{i\beta}R^\beta \right)\\
&\leq &{1\over 2}\varsigma^\delta\left(
\varsigma^{i\delta}Y_R(z;{\bf u},p,\theta)+2\varsigma ^{i\delta}c_\lambda ({\bf f})
R^\beta \right)\\
&\leq &
\varsigma^{(i+1)\delta}(Y_R(z;{\bf u},p,\theta)+c_\lambda ({\bf f})R^\beta),
\end{eqnarray*}
which completes the proof of (i) in Lemma \ref{ite}.

(ii)
Let $\rho \in ]0, \varsigma R]$ be arbitrary and $m\in\mathbb N$ be such that
\[\varsigma^{m+1}<{\rho\over R}\leq \varsigma^{m}.\]
Proceeding as in \cite{LS} we have
\begin{eqnarray*}
Y_\rho(z; {\bf u},p, \theta)& \leq&C(\varsigma)Y_{\varsigma ^{m}R}(z; {\bf u},p, \theta) 
\le (\mbox{see }(\ref{ii})\ )\\
&\le &C(\varsigma)\left({
\varsigma^{m+1}\over\varsigma}\right)^\delta(Y_R(z;{\bf u},p,\theta)+c_\lambda ({\bf f})
R^\beta)\\
&\leq &
C(\varsigma)\left({1\over \varsigma}
{\rho\over R}\right)^\delta(Y_R(z;{\bf u},p,\theta)+c_\lambda({\bf f}) R^\beta),
\end{eqnarray*}
which concludes the proof of  Lemma \ref{ite}.

\section{Proof of Theorem \ref{main}}
\label{r4}

Choose $\varsigma$ such that (\ref{tau}) is fulfilled. 
For instance, taking $\delta=\alpha/2$ and $\beta=\lambda/2$ it follows
\[\varsigma^{\alpha/2}\leq 1/C_*\quad\mbox{ and } \quad
\varsigma^{\alpha/2}+\varsigma^{\lambda/2}\leq 1.
\]
Thanks to Theorem \ref{T1_1} with the above chosen $\varsigma$ and $\beta$,
there exist $\bar\varepsilon=\varepsilon(C_*,\alpha,\lambda)$
and $\bar\gamma=\gamma(C_*,\alpha,\lambda)$ such
that for any proper weak solution $({\bf u}, p, \theta)$
of  the N-S-F problem in $Q_T$
satisfying
\begin{equation}\label{dh1}
\left.\begin{array}c
Q(z,R)\subset\subset Q_T,\quad 0<R<\bar\varepsilon,\\
R|({\bf u})_{Q(z,R)}|<1,\quad R|( {\theta})_{Q(z,R)}|<\bar\gamma/4\\
 Y_R(z; {\bf u},p, \theta)+c_\lambda ({\bf f})R^{\lambda/2} <\bar \gamma
\end{array}\right\}
\end{equation}
 and Lemma \ref{ite} can be applied.
Indeed, for all $0<\rho\leq\varsigma R$ we get
\begin{equation}\label{dh2}
Y_\rho(z;{\bf u},p,\theta)\leq C\left({\rho\over R}\right)^{\alpha/2}
( Y_R(z;{\bf u},p,\theta)+c_\lambda({\bf f}) R^{\lambda/2}).
\end{equation}

Let $\Lambda=\bar\varepsilon\bar\gamma /8$ be the desired constant.
By hypothesis (\ref{limsup}) for each $z_0\in Q_T$ 
there exists a constant $R_0>0$ that can be chosen such that $R_0<\bar\varepsilon/2$
and for all $R<R_0$ we have
\[
Q(z_0,R)\subset\subset Q_T; \quad 
R(\bar Y_R(z_0;{\bf u},p,\theta)+c_\lambda ({\bf f})R^{\lambda/2}
)<{\bar\varepsilon\bar\gamma\over 8}+{\bar\varepsilon\bar\gamma\over 8}={
\bar\varepsilon\bar\gamma\over 4}.
\]

By the continuity of $z\mapsto\bar Y_R(z;{\bf u},p,\theta) $ at $z_0$,
 there exists a neighbourhood of $z_0$, 
 ${\mathcal O}(z_0)$, such that for all $z\in {\mathcal O}(z_0)$ we have
\begin{equation}\label{gamma}
R(\bar Y_R(z;{\bf u},p,\theta)+c_\lambda({\bf f})
 R^{\lambda/2})<\bar\varepsilon\bar\gamma/4.
\end{equation}
Indeed, there exists $0<R_1\leq R_0$ such that for all 
$R<R_1$ and $Q(z,R)\subset \mathcal{O}(z_0)$
the relation  (\ref{gamma}) is satisfied.
Applying Lemma \ref{T2} with $e=(y,s)$ and $r=\bar\varepsilon$ in  (\ref{IF})-(\ref{cr})
and using  (\ref{gamma})
 we obtain
\begin{eqnarray*}
\bar Y_r(e;{\bf u}^R,p^R,\theta^R)+c_\lambda({\bf f}^R;
Q_r)r^{\lambda/2}\leq \\
\leq{R\over r}\left(
\bar Y_R(z;{\bf u},p,\theta)+
\left({R\over r}\right)^{\lambda/2}R^{\lambda/2}c_\lambda({\bf f})\right)
<\bar\gamma/4 \qquad (R<r),
\end{eqnarray*}
with $Q_r$ denoting the transported domain.
Let us consider the transported solution $({\bf u}^R,p^R,\theta^R)$
and $0<r<\min\{\bar\varepsilon,4/\bar\gamma\}$.
Then the assumption (\ref{dh1}) is fulfilled 
regarding that  $r<\bar\varepsilon$, and
\begin{eqnarray*}
r|({\bf u}^R)_{Q(e,r)}|\leq r \bar Y_r(e;{\bf u}^R,p^R,\theta^R)<1
,\\ r|( {\theta})^R_{Q(e,r)}|\leq
\bar Y_r(e;{\bf u}^R,p^R,\theta^R)<\bar\gamma/4,\\
 Y_r(e; {\bf u}^R,p^R, \theta^R)+c_\lambda({\bf f}^R;Q_r)r^{\lambda/2}\leq
4\bar Y_r(e; {\bf u}^R,p^R, \theta^R)+c_\lambda({\bf f}^R;Q_r)r^{\lambda/2} <\bar \gamma.
\end{eqnarray*}

Now using (\ref{dh2}) it results
\begin{eqnarray*}
 Y_\rho(e; {\bf u}^R,p^R, \theta^R)&\leq&
C\left({\rho\over r}\right)^{\alpha/2}\left(
 Y_r(e; {\bf u}^R,p^R, \theta^R)+c_\lambda({\bf f}^R;Q_r)r^{\lambda/2}\right)\\
&<&\bar\gamma C\left({\rho\over r}\right)^{\alpha/2}.
\end{eqnarray*}

Thus  ${\bf u}$ is H\"older continuous with exponent $\alpha/4$ in a neighbourhood
of $e$,
taking into account the parabolic version of the Campanato
criterion \cite[Theorem I.2]{camp}. 
Theorem \ref{main} is proved.

\section{Proof of Theorem \ref{dim}}
\label{fim}

We  begin by recalling the following result 
which plays a central role in the proof.
\begin{lemma}[{\cite[Lemma 11]{gg}}]\label{lgg}
Let $f\in L^1_{\rm loc}(Q_T)$ and, for $0<d<n+2$, let
 \[
F=\{z\in Q_T:\ \limsup_{R \to 0^+}R^{-d}\int_{Q (z,R)}|f|dxdt>0\}.
\]
Then, we have $\mathcal H^d(F)=0$.
\end{lemma}

Let $z\in S$  be arbitrary and, for $R>0$, denote
\[
A= {1\over\omega_n R^n}\int
 _{Q(z, R)} |{\bf u}|^{a},
\quad B={1\over\omega_n R^n}\int
 _{Q(z, R) }|p|^ {n+2\over n},\quad D={1\over\omega_n R^n}\int
 _{Q(z, R) } |\theta  |^{a'},
 \]
 with $\omega_n$ denoting the measure of the unit $n$-dimensional ball.

First, we observe that
\[S\subset \bigcup_{i=1}^5F_i,\]
where each $F_i$ is the set of points of $S$ that verify the following $i$-case.
\begin{enumerate}
\item ${[A\leq 1,\ B\leq 1,\ D\leq 1]}$
This case is impossible due to 
\begin{eqnarray*}
0<\Lambda&<& \limsup_{R \to 0^+}R(R^{-2/a}+R(R^{-2n/(n+2)}+R^{-2/a'}))\\
&=&\limsup_{R \to 0^+}(R^{1-2/a}+R^{4/(n+2)}+R^{2/a})=0,
\end{eqnarray*}
taking $a>2$ into account (cf. (\ref{defa})).
Then we find that $F_1=\emptyset$.
\item ${[A\leq 1,\ B\leq 1,\ D\geq 1]}$
Since $A^{a'/a}\leq 1$ and $B^{a'n/(n+2)}\leq 1$, we get
\begin{eqnarray*}
0<\Lambda^{a'}&< &4
\limsup_{R \to 0^+}\left(
R^{(1-2/a)a'}+R^{4a'/(n+2)}+R^{2a'/a}D\right)\\
&\leq& 4\limsup_{R \to 0^+}{R^{2a'/a}\over \omega_n R^{n}}
\int _{Q(z, R) } |\theta  |^{a'}dxdt.
\end{eqnarray*}
Then from Lemma \ref{lgg} it follows  $\mathcal H^{n-2/(a-1)}(F_2)=0$.
\item ${[A\leq 1,\ B\geq 1,\ \forall D]}$
Considering that
\[1<{2(n+2)\over n+4}<a'\leq {3\over 2}<{n+2\over n}<3\leq a<{2(n+2)\over n}\]
implies that $B^{a'n/(n+2)}\leq B$, we get
\begin{eqnarray*}
0<\Lambda^{a'}&<&4 \limsup_{R \to 0^+}\left(
R^{(1-2/a)a'}+R^{4a'/(n+2)}B+R^{2a'/a}D\right)\\
&\leq &4\limsup_{R \to 0^+}{R^{2a'/a}\over \omega_n R^{n}}
\int _{Q(z, R) } |p|^{{n+2\over n}}+|\theta  |^{a'}dxdt.\end{eqnarray*}
Then from Lemma \ref{lgg} it follows  $\mathcal H^{n-2/(a-1)}(F_3)=0$.

\item ${[A\geq 1,\ B\leq 1,\ \forall D]}$
Arguing as in the above  two cases, now
considering that $A^{a'/a}\leq A$, we get
\begin{eqnarray*}
0<\Lambda^{a'}&<& 4\limsup_{R \to 0^+}\left(
R^{(1-2/a)a'}A+R^{4a'/(n+2)}+R^{2a'/a}D\right)\\
&\leq& 4\limsup_{R \to 0^+}{R^{(1-2/a)a'}\over \omega_n R^{n}}
\int _{Q(z, R) } |{\bf u}|^a+|\theta|^{a'}dxdt.
\end{eqnarray*}
Then from Lemma \ref{lgg} it follows  $\mathcal H^{n-(a-2)/(a-1)}(F_4)=0$.

\item ${[A\geq 1,\ B\geq 1,\ \forall D]}$
Analogously to the above cases, 
now with  $A^{a'/a}\leq A$ and $B^{a'n/(n+2)}\leq B$, we get
\begin{eqnarray*}
0<\Lambda^{a'}&< &4\limsup_{R \to 0^+}\left(
R^{(1-2/a)a'}A+R^{4a'/(n+2)}B+R^{2a'/a}D\right)\\
&\leq& 4\limsup_{R \to 0^+}{R^{(1-2/a)a'}\over \omega_n R^{n}}
\int _{Q(z, R) }\left( |{\bf u}|^a+
|p|^{{n+2\over n}}+|\theta|^{a'}\right)dxdt.
\end{eqnarray*}
Then from Lemma \ref{lgg} it follows  $\mathcal H^{n-(a-2)/(a-1)}(F_5)=0$.
\end{enumerate}
Next, from the properties of measure, we have
\[\mathcal{H}^d(S)\leq\sum_{i=2}^5\mathcal{H}^d(F_i)
\leq\sum_{i=2}^5\mathcal{H}^{d_i}(F_i),
\]
supposing that $d=\max\{d_i,\ i=2,\cdots,5\}$.

Then, by the assumption (\ref{defa})
we have $a\leq 4$ $(n=2,3)$ and we conclude that $\mathcal H^{n-(a-2)/(a-1)}(S)=0$.
Therefore by definition of parabolic
Hausdorff dimension the proof of Theorem \ref{dim} is finished.

\end{document}